\documentstyle{book}
\def\Le{\Leftrightarrow}
\def\p1{\hspace{1pt}}

\begin{document}

\oddsidemargin 16.5mm
\evensidemargin 16.5mm

\thispagestyle{plain}

\noindent {\small \sc Univ. Beograd. Publ. Elektrotehn. Fak.}

\noindent {\scriptsize Ser. Mat. 8 (1997), 90--92.}

\vspace{5cc}

\begin{center}

{\Large \bf AN APPLICATION OF \boldmath $\lambda$-METHOD ON
SHAFER-FINK'S INEQUALITY \rule{0mm}{6mm}}

\vspace{1cc}

{\large \it Branko J. Male\v sevi\' c}

\vspace{1cc}

\parbox{23 cc}{{\scriptsize \bf In the paper $\lambda$-method
 Mitrinovi\' c-Vasi\' c is applied aiming to improve
Fink's inequality, and Shafer's inequality for arcus
sinus function is observed.}}
\end{center}

\vspace{1cc}

     In monography [{\bf 1}, p. 247] {\sc Shafer}'s inequality is stated:
$$ \frac{3x}{2+\sqrt{1-x^2}} \le \mbox{asin}\,x\qquad (0 \le x \le 1).
\leqno(1) $$
The equality holds only for $x=0$.
In paper [{\bf 2}] {\sc Fink} has proved the inequality:
$$ \mbox{asin}\, x \le \frac{\pi x}{2 + \sqrt{1-x^2}} \qquad
(0 \leq x \leq 1). \leqno(2) $$
The equality holds at both ends of the interval $x=0$ and $x=1$.
Let us notice that from the inequality (1) and (2) the function
$g(x) = \mbox{asin}\, x$ is bounded by the corresponding functions
from the two-parameters family of functions:
$$ \Phi_{a,b}(x) = \frac{ax}{b + \sqrt{1 - x^2}} \qquad (0\leq x\leq 1),
\leqno(3) $$
for some values of parameters $a,b > 0$. For the values of parameters
$a,b > 0$ the family $\Phi_{a,b}(x)$ is the family of raising convex functions
on variable $x$ on interval $(0,1)$. Let us apply $\lambda$-method
{\sc Mitrinovi\' c}-{\sc Vasi\' c} [{\bf 1}] on considered two-parameters family
$\Phi_{a,b}$ in order to determine the bound of function $g(x)$ under the
following conditions:
$$ \Phi_{a,b}(0)=g(0) \quad \mbox{and} \quad
\frac{{\rm d}\;}{{\rm d} x}\, \Phi_{a,b}(0) =
\frac{{\rm d}\;}{{\rm d} x}\, g(0). \leqno(4) $$
It follows that $a=b+1$. In that way we get one-parameter subfamily:
$$ f_{b}(x)=\Phi_{b+1,b}(x)=\frac{(b+1)x}{b+\sqrt{1-x^2}}
\qquad (0 \leq x \leq 1),
\leqno(5) $$
according to parameter $b>0$, which fulfills condition (4). For family
(5) the following equivalence is true:
$$ f_{b_1}(x) < f_{b_2}(x) \; \Le \; b_1 > b_2.
\leqno(6) $$

Let us consider one-parameter function of the distance:
$$ h_b(x) = f_b(x) - g(x) = \frac{(b+1)x}{b+\sqrt{1-x^2}}-\mbox{asin}\, x
\qquad (0 \le x \le 1), \leqno(7) $$
as well as its derived function:
$$\frac{{\rm d}\;}{{\rm d}x}\, h_b(x) =
\frac{\big(\sqrt{1-x^2} - (b^2-b-1)\big) \cdot x^2}{
(b+\sqrt{1-x^2}\p1 )^2 (1-x^2+\sqrt{1-x^2}\p1 )} \qquad (0 \leq x < 1).$$
Then, it holds that $h_b(0)=0$ and $\displaystyle
\frac{{\rm d}\;}{{\rm d} x} h_b(0) = 0$. In further
consideration let us use the equivalence:
$$ \frac{{\rm d}\;}{{\rm d} x} h_b(x) \geq 0 \Le \sqrt{1-x^2} \ge (b^2-b-1).
\leqno(8) $$

     The least upper bound of the function $g(x)$ from family (5), on the
basis of equivalence (6), we get for the maximum value of parameter $b$ for
which $g(x) < f_b(x)$ is true. Let us notice $\displaystyle
h_b(1) = 1 + \frac{1}{b} -\frac{\pi}{2} \ge 0$ iff $\displaystyle
b \in \Big(0,\frac{2}{\pi-2}\Big]$. If
$\displaystyle b\in \Big(0,\frac{1+\sqrt{5}}{2}\Big]$ the right side
of equivalence (8)
is always true. If $\displaystyle b \in
\Big(\frac{1+\sqrt{5}}{2},\frac{2}{\pi-2}\Big]$
the right side of equivalence (8) is true for $x \in (0,d\p1 ]$ where
$d = d(b) = \sqrt{-b^4+2b^3+b^2-2b}$. For the maximum value of parameter
$\displaystyle b_1 = \frac{2}{\pi-2}$ we find that
$d_1 = d(b_1) \cong 0.948 \in
(0,1)$. The function $h_{b_1}(x)$ fulfills $h_{b_1}(0) = h_{b_1}(1) = 0$
and reaches the maximum for $x = d_1$. Therefore for the value of the
parameter $\displaystyle b_1 = \frac{2}{\pi-2}$ the function $f_{b_1}(x)$
is the least upper bound of the function $g(x)$ from the family (5). Thus,
inequality is proved:
$$ \mbox{asin}\, x \leq \frac{\frac{\pi}{\pi-2}x}{\frac{2}{\pi-2}+\sqrt{1-x^2}}
\qquad (0 \le x \le 1). \leqno(9) $$
The equality holds at the both ends of the interval $x=0$ and $x=1$.
The maximum distance of the function $f_{b_1}(x)$ from the function $g(x)$
is reached for $x = d_1$ and it equals $h_{b_1}(d_1) \cong 0.013$. It
is directly verified that $f_{b_1}(x) = \Phi_{b_1+1,b_1}(x) < \Phi_{\pi,2}(x)$.
Thus, the given upper bound is better than the one shown in the
paper [{\bf 2}].

\newpage

The greatest lower bound of the function $g(x)$ from the family (5), on the
basis of equivalence (6), we get for the minimum value of parameter $b$ for
which it holds $f_b(x) < g(x)$. If $b \in (b_1,2)$ then the function $h_b(x)$
has a root on $(0,1)$. If $b \geq 2$ on the basis of equivalence (8) we can
conclude: $\displaystyle\frac{{\rm d}\;}{{\rm d} x}h_b(x) \le 0$. Thus, for
the value of parameter
$b_2=2$ {\sc Shafer}'s function $f_{b_2}(x)$ is the greatest lower bound
of the function $g(x)$ from the family (5).

\vspace{1cc}

\begin{center}
{\small \bf REFERENCES:}
\end{center}

\newcounter{ref}
\begin{list}{\small \arabic{ref}.}{\usecounter{ref} \leftmargin 4mm
\itemsep -1mm}

\item {\small {\sc D. S. Mitrinovi\' c, P. M. Vasi\' c:}
{\it Analytic inequalities}. Springer--Verlag 1970.}

\item {\small {\sc A. M. Fink:} {\it Two inequalities.}
Publikacije ETF Ser. Mat. {\bf 6} (1995), 48--49.}
\end{list}

{\small \hfill \break
\noindent Faculty of Electrical Engineering, \hfill
(Received May 3, 1997)\break
\noindent University of Belgrade, \hfill\break
\noindent P.O.Box 816, 11001 Belgrade, \hfill\break
\noindent Yugoslavia \hfill}

\end{document}